\documentclass[a4paper,11pt]{amsart}  

\usepackage{dacxjo} 

\title[Symmetries and equations of smooth quartics with many lines]{Symmetries and equations of smooth quartic surfaces with many lines}
\author{Davide Cesare Veniani}

\address{Institut für Mathematik \\
FB 08 - Physik, Mathematik und Informatik \\
Johannes Gutenberg-Universität \\
Staudingerweg 9, 4. OG \\
55128 Mainz, Germany}
\email{veniani@uni-mainz.de}

\date{\today}

\subjclass[2010]{%
14J28, 
14N25
}

\begin{document}

\begin{abstract}
We provide explicit equations of some smooth complex quartic surfaces with many lines, including all 10 quartics with more than 52 lines. We study the relation between linear automorphisms and some configurations of lines such as twin lines and special lines. We answer a question by Oguiso on a determinantal presentation of the Fermat quartic surface.
\end{abstract}

\maketitle

\section{Introduction}

\subsection{Principal results} In this paper, all algebraic varieties are defined over the field of complex numbers~$\CC$. If $X \subset \PP^3$ is an algebraic surface, we denote the number of lines on $X$ by~$\Phi(X)$.

While it is classically known that a smooth cubic surface contains exactly 27 lines, the number of lines on a smooth quartic surface $X$ is finite, but depends on $X$. If $X$ is general, then $\Phi(X) = 0$. Segre first claimed in 1943 that $\Phi(X) \leq 64$ for any surface~\cite{segre}. His arguments, though, contained a flaw which was corrected about 70 years later by Rams and Schütt~\cite{64lines}. Their result was then strongly improved by Degtyarev, Itenberg and Sertöz~\cite{degtyarev-itenberg-sertoz}, who showed that if $\Phi(X) > 52$, then $X$ is projectively equivalent to exactly one of a list of 10 surfaces---called $X_{64}$, $X'_{60}$, $X''_{60}$, $\bar {X}''_{60}$, $X_{56}$, $\bar {X}_{56}$, $Y_{56}$, $Q_{56}$, $X_{54}$ and $Q_{54}$---whose Néron-Severi lattices are explicitly known. The subscripts denote the number of lines that they contain. The pairs $(X_{60}'',\bar X_{60}'')$ and $(X_{56},\bar X_{56})$ are complex conjugate to each other, so these 10 surfaces correspond to 8 different line configurations.

The aim of this paper is to find an explicit defining equation of each of these 10 surfaces.

Seven of these 10 surfaces are already known in the literature. The surface~$X_{64}$---which as an immediate corollary of this list is the only surface up to projective equivalence which contains $64$ lines---was found by Schur in 1882~\cite{schur} and is given by the equation
\[
x_0(x_0^3-x_1^3) = x_2 (x_2^3 - x_3^3).
\]

The surfaces with 60 lines were found by Rams and Schütt. An equation of $X'_{60}$ is contained in~\cite{64lines}. It was found while studying a particular 6-dimensional family $\mathcal Z$ of quartics containing a line intersecting 18 or more other lines (see §\ref{subsec:symmetries-order-3}).

The surfaces $X''_{60}$ and $\bar X''_{60}$ were found by using positive characteristic methods \cite{rams-schuett-char2}. These surfaces are still smooth and contain 60 lines when reduced modulo~$2$. According to Degtyarev~\cite{degtyarev-supersingular}, 60 lines is the maximal number of lines attainable by a smooth quartic surface defined over a field of characteristic~$2$. 

The surface $X_{56}$ was studied by Shimada and Shioda~\cite{shimada-shioda} due to a peculiar property: it is isomorphic as an abstract K3 surface to the Fermat quartic surface $X_{48}$, but it is not projectively equivalent to it, as the Fermat quartic only contains $48$ lines. They provide an explicit equation of $X_{56}$ and also an explicit isomorphism between the quartic surfaces. Oguiso showed that the graph $S$ of the isomorphism $X_{48}\xrightarrow{\sim} X_{56}$ is the complete intersection of four hypersurfaces of bi-degree $(1,1)$ in $\PP^3\times \PP^3$~\cite{oguiso17}. He also asked for explicit equations of the graph $S$, which we provide in~§\ref{sec:oguiso-pairs}. As a byproduct, we obtain a determinantal description of the Fermat quartic surface.

A defining equation of the surface $Y_{56}$ is contained in~\cite{degtyarev-itenberg-sertoz}. This surface is a real surface with $56$ real lines, i.e., it attains the maximal number of real lines that can be contained in a smooth real quartic.

The three remaining surfaces are $Q_{56}$, $X_{54}$ and $Q_{54}$. In order to provide an explicit equation of $Q_{56}$ and $X_{54}$ (see §\ref{subsec:Q56} and §\ref{subsec:X54}), we investigate the group of linear automorphisms of each quartic surface, meaning the automorphisms that are restrictions of automorphisms of $\PP^3$. We use the word `symmetry' for such an automorphism. Thanks to the global Torelli theorem, the group of symmetries can be computed from the Néron-Severi group using Nikulin's theory of lattices~\cite{nikulin}. 


As for $Q_{54}$, we first find by the same method an explicit equation of $X_{52}''$, the unique surface up to projective equivalence containing configuration~$\bX_{52}''$. This surface is isomorphic to $Q_{54}$ as abstract K3~surface. Following~\cite{shimada-shioda}, we find an explicit isomorphism between the two surfaces, which in turn enables us to find an equation of $Q_{54}$ (see §\ref{subsec:Q54}).

In previous papers \cite{64lines, veniani1}, two particular configurations of lines played an important role, viz. twin lines and special lines. Their geometry is related to torsion sections of some elliptic fibrations. Here we relate these configurations to the presence of certain symmetries of order 2 and 3, thus providing yet another characterization of both phenomena (see Propositions~\ref{prop:twins-antisympl} and~\ref{prop:special-aut}).

We extend our lattice computations to all rigid (i.e., of rank~$20$) line configurations that can be found in~\cite{degtyarev-itenberg-sertoz} and~\cite{degtyarev16}. In particular, we determine for each of them the size of their group of symmetries, listed in Table~\ref{tab:main-table}. We are also able to find explicit equations of all smooth quartic surfaces containing configurations $\bX_{52}''$, $\bX_{52}'''$, $\bX_{52}^\mathrm v$, $\bY_{52}'$, $\bQ_{52}'''$ and $\bX_{51}$ (see §\ref{sec:<=52-lines}). We also provide a 1-dimensional family of quartic surfaces whose general member is smooth and contains the non-rigid configuration $\bZ_{52}$ (see §\ref{subsec:Z52}).

\subsection{Contents of the paper} In §\ref{sec:lattices} we introduce the basic nomenclature and notation about lattices. In §\ref{sec:fano-configurations} we establish the connection between lattices and configurations of lines on smooth quartic surfaces and their symmetries. In §\ref{sec:symmetries-and-lines} we investigate the properties of symmetries of order 2 and 3, relating them to some known configurations of lines, namely twin lines and special lines. In §\ref{sec:>52-lines} and §\ref{sec:<=52-lines} we explain how to find explicit equations of several quartic surfaces containing many lines. In §\ref{sec:oguiso-pairs}, we introduce Oguiso pairs and give explicit equations for the graph of the isomorphism $X_{48}\xrightarrow{\sim} X_{56}$.

\subsection*{Acknowledgements}
Thanks are due to Fabio Bernasconi, Alex Degtyarev, Noam Elkies, Dino Festi, Andrés Jaramillo Puentes, Ariyan Javanpeykar, Keiji Oguiso, Matthias Schütt and Ichiro Shimada for motivation and fruitful discussions.

\section{Lattices}
\label{sec:lattices}

Let $R$ be $\ZZ$, $\QQ$ or $\RR$. An \emph{$R$-lattice} is a free finitely generated $R$-module $L$ equipped with a non-degenerate symmetric bilinear form $\langle \ , \: \rangle \colon L \times L \rightarrow R$. 

Let $n$ be the rank of a $\ZZ$-lattice $L$ and choose a basis $\{e_1,\ldots,e_n\}$ of $L$. The matrix whose $(i,j)$-component is $\langle e_i, e_j \rangle$ is called a \emph{Gram matrix} of $L$. 
The determinant of this matrix does not depend on the choice of the basis. It is called the \emph{determinant} of $L$ and denoted $\det L$.

\subsection{Positive sign structures} By a well-known theorem of Sylvester, any $\RR$-lattice admits a diagonal Gram matrix whose diagonal entries are $\pm 1$. The numbers $s_+$ of $+1$ and $s_-$ of $-1$ are well defined and the pair $(s_+,s_-)$ is called the \emph{signature} of $L$. A lattice is \emph{positive definite} if $s_- = 0$, \emph{negative definite} if $s_+ = 0$, or \emph{indefinite} otherwise. It is \emph{hyperbolic} if $s_+ = 1$. 

A \emph{positive sign structure} of $L$ is the choice of a connected component of the manifold parameterizing oriented $s_+$-dimensional subspaces $\Pi$ of $L$ such that the restriction of $\langle \ , \: \rangle$ to $\Pi$ is positive definite.

By definition, the signature and the positive sign structure of a $\ZZ$- or $\QQ$-lattice $L$ are those of $L\otimes \RR$. 

\subsection{Discriminant forms} Let $D$ be a finite abelian group. A \emph{finite symmetric bilinear form} on $D$ is a a homomorphism $b\colon D \times D \rightarrow \QQ/\ZZ$ such that $b(x,y) = b(y,x)$ for any $x,y \in D$. A \emph{finite quadratic form} on $D$ is a map $q\colon D \rightarrow \QQ/2\ZZ$ such that
\begin{enumerate}[(i),noitemsep]
\item $q(nx) = n^2q(x)$ for $n \in \ZZ$, $x \in D$; 
\item the map $b\colon D\times D \rightarrow \QQ/\ZZ$ defined by $(x,y) \mapsto (q(x+y) - q(x) - q(y))/2$ is a finite symmetric bilinear form.
\end{enumerate}
A finite quadratic form $(D,q)$ is \emph{non-degenerate} if the associated finite symmetric bilinear form $b$ is non-degenerate. We denote by $\Orth(D,q)$ the group of automorphisms of a non-degenerate finite quadratic form $(D,q)$.

Let $L$ be a $\ZZ$-lattice. The \emph{dual lattice} $L^\vee$ of $L$ is the group of elements $x \in L \otimes \QQ$ such that $\langle x,v\rangle \in \ZZ$ for all $v \in L$. The dual lattice $L^\vee$ is a free $\ZZ$-module which contains $L$ as a submodule of finite index. In particular, $L$ and $L^\vee$ have the same rank.

A $\ZZ$-lattice $L$ is \emph{even} if $\langle x,x \rangle \in 2\ZZ$ holds for any $x \in L$. 
Given an even $\ZZ$-lattice $L$, the \emph{discriminant form} $(D_L,q_L)$ of $L$ is the finite quadratic form on the group $D_L := L^\vee/L$ defined by $q_L(\bar x) = \langle x,x\rangle \mod 2\ZZ$, where $\bar x \in D_L$ denotes the class of $x \in L^\vee$ modulo $L$.

\subsection{Genera}


Let $(s_+,s_-)$ be a pair of non-negative integers and $(D,q)$ be a non-degenerate finite quadratic form. The \emph{genus} $\cG$ determined by $(s_+,s_-)$ and $(D,q)$ is the set of isometry classes of even $\ZZ$-lattices $L$ of signature $\sign(L) = (s_+,s_-)$ and discriminant form $(D_L,q_L)\cong (D,q)$. 

The \emph{oriented genus} $\cG^\orient$ determined by $(s_+,s_-)$ and $(D,q)$ is the set of equivalence classes of pairs $(L,\theta)$, where $L$ is a lattice whose isometry class belongs to $\cG$, and $\theta$ is a positive sign structure on $L$. We say that $(L,\theta)$ and $(L',\theta')$ are equivalent if there is an isometry $L \xrightarrow{\sim} L'$ which maps $\theta$ to $\theta'$. There is an obvious forgetful map $\cG^\orient \rightarrow \cG$.


\subsection{Positive definite lattices of rank 2} 
\label{subsec:pos-def-rank-2}
If $T$ is a positive definite even $\ZZ$-lattice of rank 2 and discriminant form $(D,q)$, then there exists a unique triple $(a,b,c)$ of integers with $0< a \leq c$, $0 \leq 2b \leq a$, $a$ and $c$ even, such that $T$ admits a Gram matrix of the form
\[
\left[\begin{array}{cc}   
 a & b \\ 
 b & c \\
\end{array}\right].
\]
We denote by $[a,b,c]$ the isometry class of $T$.

It is easy to compute the genus $\cG$ and oriented genus $\cG^\orient$ determined by $(2,0)$ and $(D,q$), since $\frac34 c \leq \det T \leq c^2$. The preimage of $[a,b,c]$ under the map $\cG^\orient\rightarrow \cG$ has either one or two elements. It has one element if and only if $T$ admits an orientation reversing autoisometry, which is the case if and only if $a = c$, or $a = 2b$, or $b = 0$.

\section{Fano configurations}
\label{sec:fano-configurations}

A \emph{$d$-polarized lattice} is a hyperbolic lattice $S$ together with a distinguished vector $h \in S$ such that $h^2 = d$, called the \emph{polarization}. A $\emph{line}$ in a polarized lattice $(S,h)$ is a vector $v \in S$ such that $v^2 = -2$ and $v\cdot h = 1$. The set of lines is denoted by $\Fn(S,h)$.
A \emph{configuration} is a $4$-polarized lattice $(S,h)$ which is generated over $\QQ$ by $h$ and all lines in $\Fn(S,h)$. A configuration is \emph{rigid} if $\rank S = 20$.

Let $X \subset \PP^3$ be a smooth quartic surface. The primitive sublattice $\cF(X)$ of $H^2(X,\ZZ)$ spanned over $\QQ$ by the plane section $h$ and the classes of all lines on $X$ is called the \emph{Fano configuration} of $X$. The plane section defines a polarization of $\cF(X)$. A configuration is called \emph{geometric} if it is isometric as a polarized lattice to the Fano configuration of some quartic surface.

\subsection{Projective equivalence classes} 
\label{subsec:proj-equiv-classes} 
Let $S$ be a rigid geometric configuration. Consider the (non-empty) set of projective equivalence classes of smooth quartic surfaces whose Fano configuration is isometric to $S$. This set is finite and its cardinality can be computed in the following way. Let $\cG_S$ and $\cG^\orient_S$ be the (oriented) genus determined by $(2,0)$ and $(D_S,-q_S)$. Fix a positive definite lattice $T$ of rank $2$ whose class is in $\cG_S$, and let $\psi\colon (D_T,q_T) \rightarrow (D_S,-q_S)$ be an isomorphism. We can identify $\Orth_h(S)$ as a subgroup of $\Orth(D_T,q_T)$. Define 
\[ \cl(S) := |\Orth^+(T) \backslash \Orth(D_T,q_T) / \Orth_h(S)|.\]
This number does not depend on the $T$ and $\psi$ chosen. The number of projective equivalence classes is equal to 
\[\Cl(S) := \cl(S) \cdot |\cG^\orient_S|.\] 
For more details, see \cite[Remark 3.6]{degtyarev16}.

\subsection{Symmetries} \label{subsec:symmetries}

Let $(S,h)$ be geometric Fano configuration. Fix a lattice $T$ in $\cG_S$ and an isomorphism $\psi\colon (D_S,-q_S) \xrightarrow{\sim} (D_T,q_T)$. Let $\Gamma_T$ be the image of $\Orth^+(T)$ in $\Orth(D_T,q_T)$. We define the subgroup 
\[
\Gamma_S := \{ \psi^{-1}\gamma \psi \mid \gamma \in \Gamma_T\} \subset \Orth(D_S,q_S).
\]
The group $\Gamma_S$ does not depend on the $T$ and $\psi$ chosen (see \cite[§2.4]{degtyarev16}).

Let $\eta_S\colon \Orth_h(S) \rightarrow \Orth(D_S,q_S)$ be the natural homomorphism. We consider the subgroup of \emph{symmetries} of~$S$
\[ \Sym(S) := \{ \varphi \in \Orth_h(S) \mid  \eta_S(\varphi) \in \Gamma_S \} \]
and the subgroup of \emph{symplectic symmetries} of~$S$
\[
	\Sympl(S) := \{ \varphi \in \Orth_h(S) \mid \eta_S(\varphi) = \id \}.
\]

Let $X \subset \PP^3$ be a smooth quartic surface. 
 A \emph{symmetry} of $X$ is an automorphism $\varphi\colon X\rightarrow X$ which is the restriction of an automorphism of $\PP^3$. A symmetry $\varphi$ is \emph{symplectic} if it acts trivially on $H^{2,0}(X)$. We denote the group of symmetries of $X$ by $\Sym(X)$ and the subgroup of symplectic symmetries by $\Sympl(X)$.


Let $S := \cF(X)$ be the Fano configuration of $X$. A symmetry of $X$ induces a symmetry of $S$, thus giving a homomorphism $\Sym(X) \rightarrow \Sym(S)$. The following proposition is a consequence of Nikulin's theory of lattices and the global Torelli theorem.

\begin{proposition} \label{prop:torelli}
If $S = \cF(X)$ is a rigid geometric Fano configuration, then the homomorphism $\Sym(X) \rightarrow \Sym(S)$ is an isomorphism. Furthermore, symplectic automorphisms of~$X$ correspond to symplectic automorphisms of~$S$ under this isomorphism.
\end{proposition}

\subsection{The main table}
\label{subsec:main-table}
Table~\ref{tab:main-table} lists all known rigid geometric Fano configurations found in \cite{degtyarev-itenberg-sertoz} and \cite{degtyarev16}. Let $N := |\Fn(S)|$. Note that $N$ is always equal to the subscript in the name of the configuration.
\begin{itemize}
\item Since $S$ is rigid, an element in $\Orth_h(S)$ corresponds uniquely to a permutation $\sigma \in \Symmetricgroup_N$ such that $\langle l_i ,\, l_j\rangle = \langle l_{\sigma(i)},\, l_{\sigma(j)} \rangle$.
\item The sixth column contains a list of all elements of $\cG_S$ (see \ref{subsec:pos-def-rank-2}). Those classes that correspond to two elements in $\cG_S^\orient$ are marked by an asterisk~${}^*$. We write ${}^{\times 2}$ if $\cl(S) = 2$; otherwise, $\cl(S) = 1$.
\item We compute $\Sym(S)$ and $\Sympl(S)$ using the definition.
\end{itemize}

The list of rigid configurations with exactly $52$ lines is not known to be complete. The only known non-rigid configuration with $52$ lines is $\bZ_{52}$ (see~§\ref{subsec:Z52}). On the other hand, there are certainly many more configurations with less than $52$ lines than the ones listed here. 

For our computations we used \verb|GAP| \cite{GAP4}.

\begin{table}
\caption{Rigid geometric Fano configurations with many lines. For an explanation of the entries, see~§\ref{subsec:main-table}.}
\label{tab:main-table}
\begin{tabular}{ccccrlc}
\toprule
$S$ & $|\Orth_h(S)|$ & $|\Sym(S)|$ & $|\Sympl(S)|$ & $\det T$ & $T^{\times \cl(S)}$ &  see \\
\midrule 
$\bX_{64}$ & 4608 & 1152 & 192 		& 48 & $[8,4,8]$  & \cite{schur} \\ 
$\bX_{60}'$ & 480 & 120 & 60 		& 60 & $[4,2,16]$ & \cite{64lines} \\ 
$\bX_{60}''$ & 240 & 120 & 60 		& 55 & $[4,1,14]^*$ & \cite{rams-schuett-char2} \\ 
$\bX_{56}$ & 128 & 64 & 16 			& 64 & $[8,0,8]^{\times 2}$  & \cite{shimada-shioda}, §\ref{sec:oguiso-pairs} \\
$\bY_{56}$ & 64 & 32 & 16 			& 64 & $[2,0,32]$ & \cite{degtyarev-itenberg-sertoz} \\
$\bQ_{56}$ & $384$ & $96$ & $48$ 	& 60 & $[4,2,16]$ & §\ref{subsec:Q56} \\ 
$\bX_{54}$ & $384$ & $48$ & $24$ 	& 96 & $[4,0,24]$ & §\ref{subsec:X54} \\
$\bQ_{54}$ & $48$ & $8$ & $8$ 		& 76 & $[4,2,20]$ & §\ref{subsec:Q54}  \\ 
\midrule
$\bX_{52}'$ & 24 & 3 & 3 			& 80 & $[8,4,12]$ & 	\\
$\bX_{52}''$ & 36 & 6 & 6 			& 76 & $[4,2,20]$ &	§\ref{subsec:X52ii} 
\\ 
$\bX_{52}'''$ & 320 & 80 & 20 		& 100 & $[10,0,10]$ & §\ref{subsec:X52iii} \\
$\bX_{52}^\mathrm v$ & 32 & 8 & 4 	& 84 & $[10,4,10]$ & §\ref{subsec:X52v} \\ 
$\bY_{52}'$ & 8 & 8 & 4 			& 76 & $\begin{cases} [2,0,38] \\ [8,2,10]^* \end{cases}$ & §\ref{subsec:Y52i} \\ 
$\bY_{52}''$ & 8 & 8 & 4 			& 79 & $\begin{cases} [2,1,40] \\ [4,1,20]^* \\ [8,1,10]^* \end{cases}$ & \\ 
$\bQ_{52}'$ & 64 & 8 & 8 			& 96 & $[4,0,24]$ & \\
$\bQ_{52}''$ & 64 & 16 & 16 		& 80 & $[8,4,12]^{\times 2}$ & \\ 
$\bQ_{52}'''$ & 96 & 24 & 4 		& 75 & $[10,5,10]$ & §\ref{subsec:Q52iii} \\ 

\midrule

$\bX_{51}$ & 12 & 6 & 6 			& 87 & $\begin{cases} [4,1,22]^* \\ [6,3,16] \end{cases}$ & §\ref{subsec:X51} \\ 
$\bX_{50}'$ & 18 & 3 & 3 			& 75 & $[4,2,28]$ & \\
$\bX_{50}''$ & 12 & 3 & 3 			& 96 & $[4,0,24]^{\times 2}$ & \\
$\bX_{50}'''$ & 16 & 4 & 2 			& 96 & $[4,0,24]^{\times 2}$ & \\
$\bX_{48}$ & 6144 & 1536 & 384 		& 64 & $[8,0,8]$  & §\ref{sec:oguiso-pairs} \\ 
$\bY_{48}'$ & 8 & 2 & 1  			& 96 & $[2,0,48]$ & \\
$\bY_{48}''$ & 8 & 4 & 4			& 95 & $\begin{cases} [2,1,48] \\ [8,1,12]^* \\ [10,5,12] \end{cases}$ & \\
$\tilde \bY_{48}'$ & 48 & 48 & 24 	& 76 & $\begin{cases} [2,0,38] \\ [8,2,10]^* \end{cases}$ & \\
$\tilde \bY_{48}''$ & 12 & 12 & 6 	& 79 & $\begin{cases} [2,1,40] \\ [4,1,20]^* \\ [8,1,10]^* \end{cases}$ & \\
$\bQ_{48}$ & 128 & 16 & 8 			& 80 & $[8,4,12]$ & \\
\bottomrule
\end{tabular}
\end{table}

\section{Symmetries and lines on quartic surfaces}
\label{sec:symmetries-and-lines}

In this section we recall some basic facts about symplectic automorphisms of K3 surfaces. Moreover, we study symmetries of smooth quartic surfaces of order 2 and 3. Some of these symmetries are related to particular configurations of lines, which have played a major role in previous works (cf. \cite{64lines,veniani1}): twin lines and special lines.

\subsection{Symplectic automorphisms} \label{subsec:sympl-autom}

We recall some basic properties of symplectic automorphisms of K3 surfaces. For more details, see~\cite{nikulin-finite-groups}.

If $\varphi\colon Y \rightarrow Y$ is an automorphism of an algebraic variety $Y$ its \emph{fixed locus} is denoted by $\Fix(\varphi)$. If $\varphi\colon X\rightarrow X$ is a symmetry of a quartic surface $X$, then we denote by the same letter also the corresponding automorphism of~$\PP^3$. If it is not clear from the context to which fixed locus we are referring, we write $\Fix(\varphi,X)$ or $\Fix(\varphi,\PP^3)$.

Let $\varphi\colon Y\rightarrow Y$ be a symplectic automorphism of a K3 surface $Y$. If the order $n$ of $\varphi$ is finite, then $n\leq 8$ and $\Fix(\varphi)$ consists of a finite number $f_n$ of points. This number $f_n$ depends only on the order $n$. The following list shows all values of $f_n$ for $n = 1,\ldots,8$. 
\begin{table}[h]
\begin{tabular}{c|ccccccc}
$n$ & 2 & 3 & 4 & 5 & 6 & 7 & 8 \\
$f_n$ & 8 & 6 & 4 & 4 & 2 & 3 & 2 \\
\end{tabular}
\end{table}

\subsection{Type of a line} Let $l$ be a line on a smooth quartic surface $X$. Consider the pencil of planes $\{\Pi_t\}_{t\in \PP^1}$ containing $l$. The curve $C_t := \Pi_t \cap X$ is called the \emph{residual cubic} in the plane $\Pi_t$. If $C_t$ splits into three lines, then $C_t$ is called a \emph{$3$-fiber}. If $C_t$ splits into a line and an irreducible conic, then $C_t$ is called a \emph{$1$-fiber}. The \emph{type} of the original line $l$ is the pair $(p,q)$, where $p$ (resp. $q$) is the number of $3$-fibers (resp. $1$-fibers) of $l$.

The name ``fiber'' comes from the fact that the morphism $X\rightarrow \PP^1$ whose fiber over $t \in \PP^1$ is $C_t$ is an elliptic fibration. In Kodaira's notation, a $3$-fiber corresponds to a fiber of type $\I_3$ or $\IV$, while a $1$-fiber corresponds to a fiber of type $\I_2$ or $\III$. The \emph{discriminant} of a line is the discriminant of its induced elliptic fibration. It is a homogeneous polynomial of degree~24 in two variables.

The restriction of $X\rightarrow \PP^1$ to $l$ is a separable morphism of curves of degree~3. Its ramification divisor has degree~$4$. The \emph{ramification points of $l$} are the ramification points with respect to this morphism.

\subsection{Symmetries of order 2 and twin lines}
Once coordinates on $\PP^n$ are chosen, any automorphism of $\PP^n$ (hence, any symmetry of a quartic surface) is represented by a square matrix of size $n$, well defined up to multiplication by a non-zero scalar. The diagonal matrix with entries $a,b,c,d \in \CC$ will be denoted by $\diag(a,b,c,d)$. The following lemma is an easy consequence of the fact that a complex matrix of finite order is diagonalizable.

\begin{lemma} \label{lem:autP3-ord2}
If $\varphi\colon \PP^3 \rightarrow \PP^3$ is an automorphism of order $2$, then $\Fix(\varphi)$ consists of the disjoint union of either two lines, or one point and one plane.
\end{lemma}

\begin{proposition} \label{prop:fix-sympl-ord-2}
If $\sigma\colon X \rightarrow X$ is a symplectic symmetry of order~$2$, then there exist two lines $l_1$ and $l_2$ in $\PP^3$ such that $\Fix(\sigma,\PP^3) = l_1\cup l_2$. Moreover, each $l_i$ intersects $X$ in $4$ distinct points. 
\end{proposition}
\proof A symplectic automorphism of order $2$ of a K3 surface has exactly 8 fixed points, so the statement follows from Lemma~\ref{lem:autP3-ord2} (cf. \cite{vanGeemenSarti}). \endproof

Two disjoint lines $l',l''$ on $X$ are \emph{twin lines} if there exist $10$ other distinct lines on $X$ which intersects both $l'$ and $l''$.

\begin{proposition} \label{prop:twins-antisympl}
If $l',l''$ are two disjoint lines on $X$, then the following conditions are equivalent:
\begin{enumerate}[(a)]
\item $l'$ and $l''$ are twin lines;
\item there exists a non-symplectic symmetry $\tau\colon X \rightarrow X$ of order~$2$ such that $\Fix(\tau,X) = \Fix(\tau,\PP^3) = l' \cup l''$.
\end{enumerate}

\end{proposition}
\proof The implication (a) $\implies$ (b) follows from \cite[Remark 3.4]{veniani1}.

Suppose that (b) holds. Up to coordinate change we can suppose that $l',l''$ are given by $x_0 = x_1 = 0$ and $x_2 = x_3 = 0$, respectively. Then, necessarily $\tau = \diag(1,1,-1,-1)$. Imposing that $\tau$ is a non-symplectic symmetry of $X$ leads to the vanishing of nine coefficients, so $X$ belongs to the family $\mathcal A$ in \cite[Proposition 3.2]{veniani1}. By the same proposition, $l'$ and $l''$ are twin lines.
\endproof

\begin{corollary} \label{cor:p-fiber-twin-line}
Suppose that $l',l''$ are twin lines and that the residual cubic in a plane containing $l'$ splits into three lines $m_1,m_2,m_3$. If $m_1$ intersects~$l''$, then the intersection point of $m_2$ and $m_3$ lies on $l'$.
\end{corollary}
\proof The non-symplectic symmetry $\tau$ fixes $m_1$ as a set, but not pointwise. Moreover, $\tau$ exchanges $m_2$ and $m_3$, otherwise there would be at least 3 fixed points on $m_1$. As $l'$ is fixed pointwise, the intersection point of $m_2$ and $m_3$ lies on $l'$. 
\endproof

\subsection{Symmetries of order 3 and special lines} \label{subsec:symmetries-order-3}

In a similar way as for order 2, we can prove the following lemma.

\begin{lemma} \label{lem:autPn-ord3}
Let $\varphi\colon \PP^n \rightarrow \PP^n$ be an automorphism of order $3$.
\begin{itemize}[(a),noitemsep]
\item If $n = 2$, then $\Fix(\varphi)$ consists of the disjoint union of either one point and one line, or three points.
\item If $n = 3$, then $\Fix(\varphi)$ consists of the disjoint union of either one point and one plane, or two lines, or two points and one line.
\end{itemize}
\end{lemma}

A line $l$ on a smooth quartic $X$ is said to be \emph{special}, if one can choose coordinates so that $l$ is given by $x_0 = x_1 = 0$ and $X$ belongs to the following family found by Rams--Schütt~\cite{64lines}:
\begin{equation} \label{eq:familyZ}
\mathcal Z\colon x_0x_3^3 + x_1x_2^3 + x_2x_3q(x_0,x_1) + g(x_0,x_1) = 0,
\end{equation}
where $q$ and $g$ are polynomials of degree $2$ and $4$, respectively.

\begin{proposition} \label{prop:special-aut}
If $l$ is a line on a smooth quartic surface $X$, then the following conditions are equivalent:
\begin{enumerate}[(a)]
\item the line $l$ is special;
\item there exists a symplectic symmetry $\sigma\colon X \rightarrow X$ of order~$3$ which preserves each plane containing $l$ as a set.
\end{enumerate}
\end{proposition}
\proof 
If (a) holds, then $\sigma = \diag(1,1,\zeta\,\zeta^2)$, where $\zeta$ is a primitive 3rd root of unity, is the required symmetry.

Conversely, suppose that (b) holds. Since $\sigma$ is symplectic, $\Fix(\sigma,X)$ consists of $6$ distinct points. As $\sigma$ preserves $l$ as a set, $\sigma$ exactly two fixed points on $l$, say $P_1$ and $P_2$. As all ramification points of $l$ are fixed, $P_1$ and $P_2$ are the only ramification points, necessarily of ramification index~$3$. Moreover, $\Fix(\sigma,\PP^3)$ cannot be the disjoint union of one plane $\Pi$ and one point, because the curve $\Pi \cap X$ would be fixed pointwise. 

Suppose that all lines in $\Fix(\sigma,\PP^3)$ pass through $P_1$ or $P_2$. Choose a plane $\Pi$ containing $l$, but not containing any lines in $\Fix(\sigma,\PP^3)$. Then, $\Fix(\sigma,\Pi) = \{P,Q\}$, but this contradicts Lemma~\ref{lem:autPn-ord3}. 

Hence, $\Fix(\sigma,\PP^3)$ is the disjoint union of $P_1$, $P_2$ and one line $m$. For $i=1,2$, let $\Pi_i$ be the plane containing $l$ whose residual cubic contains $P_i$ for $i=1,2$, and let $Q_i$ be the point of intersection of $\Pi_i$ with $m$ for $i=1,2$. If we choose coordinates so that $P_1 = (0,0,1,0)$, $P_2 = (0,0,0,1)$, $Q_1 = (1,0,0,0)$ and $Q_2 = (0,1,0,0)$, then $\sigma = \diag(1,1,\zeta\,\zeta^2)$. In particular, $l$ is the line $x_0 = x_1 = 0$. By an explicit computation, imposing that $\sigma$ is a symplectic symmetry of $X$, we obtain that $X$ is a member of family $\mathcal Z$, i.e., $l$ is special.
\endproof

The richness of the geometry of special lines comes from the presence of a $3$-torsion section on a certain base change of the fibration induced by $l$. We will use the following fact.

\begin{proposition}[\cite{64lines}] \label{prop:(6,q)-is-special}
If $l$ is a line of type $(6,q)$, with $q > 0$, then $l$ is special. Moreover, each $1$-fiber intersects $l$ at a single point.
\end{proposition}

\section{Configurations with more than 52 lines}
\label{sec:>52-lines}

In this section we compute equations for $Q_{56}$, $X_{54}$ and $Q_{54}$. For our calculations we used \verb|SageMath|~\cite{sagemath} and \verb|Singular|~\cite{DGPS}. Discriminants of elliptic fibrations induced by lines are computed with the formulas found in~\cite{jac-of-gen-1-curves}.

\subsection{Conventions} \label{subsec:conventions}
As a general approach, given a rigid geometric configuration $S$, we let $X$ be a smooth quartic surface containing $S$, i.e., such that $\cF(X)$ is isometric to $S$. The complete knowledge of $\Sym(X)$ and $\Sympl(X)$ comes from the complete knowledge of $\Sym(S)$ and $\Sympl(S)$ thanks to Proposition~\ref{prop:torelli}.

A \emph{quadrangle} is a set of four non-coplanar lines $l_0,\ldots,l_3$ such that $l_i$ intersects $l_{i+1}$ for each $i=0,\ldots,3$ (with subscripts interpreted modulo~$4$). The points of intersection of $l_i$ and $l_{i+1}$ are called \emph{vertices} of the quadrangle. A \emph{star} is a set of four coplanar lines intersecting in one single point.
When a symmetry $\varphi$ of a quartic surface $X$ fixes a line, a quadrangle, etc. on $X$, it is understood that $\varphi$ fixes the line, the quadrangle, etc. as a set, unless we explicitly state that $\varphi$ fixes them pointwise.

%


\subsection{Configuration \texorpdfstring{$\bQ_{56}$}{Q56}} \label{subsec:Q56}

Let $Q_{56}$ be a quartic surface containing configuration~$\bQ_{56}$. Then, there exist four non-symplectic symmetries of order~$2$ $\tau_0,\ldots,\tau_3$ such that $\tau_i\tau_j = \tau_i\tau_j$ for all $i,j$ (there are three such quadruples). Since $\tau_i$ does not fix two lines on $X$ and the fixed locus of a non-symplectic automorphism of order~$2$ on a K3~surface does not contain isolated points, $\Fix(\tau_i,\PP^3)$ consists of a point $P_i$ (not belonging to $X$) and a plane $\Pi_i$. As $\tau_i$ and $\tau_j$ commute, $\tau_j(P_i) = P_i$ and $\tau_j(\Pi_i) = \Pi_i$. For $j\neq i$, $\tau_j$ does not fix $\Pi_i$ pointwise (otherwise it would coincide with $\tau_i$), so $\Fix(\tau_j,\Pi_i)$ is the union of a point and a line. It follows that $P_j \in \Pi_i$ and that the four planes $\Pi_0,\ldots,\Pi_3$ are in general position. Up to coordinate change, we can suppose that $\Pi_i$ is given by $x_i = 0$. Thus, $\tau_0 =\diag(-1,1,1,1),\ldots,\tau_3 =\diag(1,1,1,-1)$. 

Up to relabeling, both $\sigma_1 := \tau_0\tau_2 = \tau_1\tau_3$ and $\sigma_2:=\tau_0\tau_3 = \tau_1\tau_2$ fix eight lines on $Q_{56}$, say $a_0,\ldots,a_7$ and $b_0,\ldots,b_7$, respectively. The lines $a_0,\ldots,a_7$ form two quadrangles $\{a_0,\ldots,a_3\}$ and $\{a_4,\ldots,a_7\}$. Hence, the vertices of these quadrangles are the eight fixed points of $\sigma_1$ on $X$. By Proposition~\ref{prop:fix-sympl-ord-2}, they lie on the lines $l\colon x_0 = x_2 = 0$ and $l'\colon x_1 = x_3 = 0$. Analogously, the lines $b_0,\ldots,b_7$ form two quadrangles $\{b_0,\ldots,b_3\}$ and $\{b_4,\ldots,b_7\}$, whose vertices lie on $m\colon x_0 = x_3 = 0$ and $m'\colon x_1 = x_2 = 0$. 

By inspection of $\bQ_{56}$ and its symmetry group, we see that (up to relabeling)
\begin{itemize}
\item $a_i$ intersects $b_j$ if and only if $i \equiv j \mod 2$;
\item there exists a symplectic automorphism $\sigma'_1$ of order~$2$ which fixes $a_0,a_2,a_4,a_6$;
\item there exists a symplectic automorphism $\sigma'_2$ of order~$2$ which fixes $b_0,b_2,b_4,b_6$.
\end{itemize}

Up to rescaling variables and coefficients, we have
\[ \sigma'_1 = \left[\begin{array}{cccc}   
   &   &   & 1 \\ 
   &   & 1 &   \\
   & 1 &   &   \\
 1 &   &   &    \\ 
\end{array}\right]\qquad \text{and} \qquad
\sigma'_2 = \left[\begin{array}{cccc}   
  &   & 1 &   \\ 
  &   &   & 1 \\
1 &   &   &   \\
  & 1 &   &   \\ 
\end{array}\right]. \]

If $a_0$, $a_1$ intersect each other at the point $(0,p,0,1)$, and $b_0$, $b_1$ intersect each other at the point $(0,q,1,0)$, for $p,q\in \CC$, then the surface belongs to the $2$-dimensional family given by
\begin{multline*}
2 \, p^{2} q^{2} {\left(x_{0}^{4} + x_{1}^{4} + x_{2}^{4} + x_{3}^{4}\right)}  + {\left(p^{4} + 1\right)} {\left(q^{4} + 1\right)} {\left(x_{0}^{2} x_{1}^{2} + x_{2}^{2} x_{3}^{2}\right)} \\ 
- 2 \, q^{2}{\left(p^{4} + 1\right)} {\left(x_{0}^{2} x_{2}^{2} + x_{1}^{2} x_{3}^{2}\right)} - 2 \,p^{2} {\left(q^{4} + 1\right)} {\left(x_{0}^{2} x_{3}^{2} + x_{1}^{2} x_{2}^{2}\right)} = 0.
\end{multline*}

The automorphism $\sigma'_1$ fixes also two lines $c,d$ which form a quadrangle with $a_0$ and $a_4$; necessarily, $c$ and $d$ pass through the fixed points of $\sigma'_1$ on $a_0$ and $a_4$, which can be explicitly computed. Up to exchanging $p$ with $-p$ or $q$ with $-q$ (which does not influence the equation of the surface), it follows that
\[q  = \frac{p-1}{p+1}. \]

Finally, the residual conic in the plane containing $a_0$ and $b_0$ is reducible. This means that $p$ satisfies
\[
p^4 - p^3 + 2\,p^2 + p + 1 = 0.
\]

\begin{remark} The surface $Q_{56}$ itself is defined over $\QQ(\sqrt{-15})$. All lines are defined over $\QQ(p)$. The surface contains 24 lines of type $(3,7)$, whose fibrations have one singular fiber of type $\mathrm{III}$.


\end{remark}

\subsection{Configuration \texorpdfstring{$\bX_{54}$}{X54}}
\label{subsec:X54}

Let $X_{54}$ be a quartic surface whose Fano configuration is isometric to $\bX_{54}$. Then $X_{54}$ contains 4 special lines of type $(6,2)$ and 10 pairs of twin lines. 
In particular, there is a quadrangle containing two opposite lines of type $(6,2)$, say $l_0$ and $l_2$, and a pair of twin lines of type $(0,10)$, say $l_1$ and $l_3$. 
Up to coordinate change, we can suppose that $l_i$ is the line $x_i = x_{i+1} = 0$. The non-symplectic symmetry $\sigma_1$ corresponding to the twin pair formed by $l_1$ and $l_3$ is $\sigma_1 = \diag(1,-1,-1,1)$ (see Proposition~\ref{prop:twins-antisympl}).

By inspection of $\Sym(\bX_{54})$, we find two symplectic symmetries $\sigma_2$ and $\sigma_3$ with the following properties
\begin{itemize}
\item $\sigma_2(l_0) = l_0,\,\sigma_2(l_1) = l_3,\,\sigma_2(l_2) = l_2$;
\item $\sigma_3(l_0) = l_2,\,\sigma_3(l_1) = l_1,\, \sigma_3(l_3) = l_3$.
\end{itemize}
Therefore, there exist $a,b,c,d\in\CC$ such that
\begin{equation*} \label{eq:symmetric-quadrangle}
\sigma_2 = \left[\begin{array}{cccc}   
   &   1 &   & \\
ab &&& \\
 &&&   a \\
 &&   b   & \\
\end{array}\right], \qquad 
\sigma_3 = \left[\begin{array}{cccc}   
 &  &   & 1 \\ 
 &  & c &   \\
 & d &  & \\
cd & & & 
\end{array}\right].
\end{equation*}

By Proposition~\ref{prop:(6,q)-is-special}, $l_0$ is a special line. This implies that $c = -b$ and that the residual conic is tangent to $l_0$ at the point of intersection with $l_1$ (see Proposition~\ref{prop:(6,q)-is-special}). 

Imposing all these conditions and normalizing the remaining coefficients, we find that $X_{54}$ is defined by the following equation:
\begin{multline*}
3 \, x_{0}^{3} x_{2} - 3 \, x_{0} x_{1} x_{2}^{2} - x_{0} x_{2}^{3} + 3 \, x_{1}^{3} x_{3}  + 3 \, x_{0}^{2} x_{2} x_{3} + 3 \, x_{1}^{2} x_{2} x_{3} - 3 \, x_{0} x_{1} x_{3}^{2} - x_{1} x_{3}^{3} = 0.
\end{multline*}

If $\xi$ is a primitive 12th root of unity, the lines in the plane $x_0 = \xi^3x_1$ are all defined over $\QQ(\xi)$. One can check that the three lines other than $l_0$ in this plane have type $(2,8)$.

\subsection{Configuration \texorpdfstring{$\bQ_{54}$}{Q54}} \label{subsec:Q54}
Let $X_{54}$ be a surface containing configuration~$\bQ_{54}$. All symmetries of $X_{54}$ are symplectic of order~$2$. There is only one symmetry~$\sigma$ which fixes $4$ disjoint lines $l_0,\ldots,l_3$. Observe that the restriction of~$\sigma$ to each $l_i$ must have 2 fixed points: these are all the 8 fixed points of $\sigma$ on~$X_{54}$. By Proposition~\ref{prop:fix-sympl-ord-2}, these $8$ points lie on two lines $l',l''$ in $\PP^3$. 

There are then $3$ more symmetries $\tau_1,\tau_2,\tau_3$, which fix two lines $m_1,m_2$ on $X_{54}$ and act in this way on $l_0,\ldots,l_3$:
\begin{itemize}
\item $\tau_1(l_0) = l_1,\,\tau_1(l_2) = l_3$;
\item $\tau_2(l_0) = l_2,\,\tau_2(l_1) = l_3$.
\item $\tau_3(l_0) = l_3,\,\tau_3(l_1) = l_2$.
\end{itemize}

Since $\Sym(\bQ_{54})$ is a commutative group, each $\tau_i$ permutes $l'$ and $l''$. As $\tau_3 = \tau_1 \circ \tau_2$, at least one $\tau_i$ fixes $l'$ and $l''$; hence, by symmetry, each $\tau_i$ fixes $l'$ and $l''$. 
Let $m', m''$ be the lines fixed pointwise in $\PP^3$ by $\tau_1$. The lines $l',m',l'',m''$ form a quadrangle. Up to coordinate change, we can suppose that
$l'\colon x_1 = x_3 = 0, l''\colon x_0 = x_2 = 0, m'\colon x_0 = x_1 = 0, m''\colon x_2 = x_3 = 0$.

In these coordinates, $\sigma = \diag(1,-1,1,-1)$ and $\tau_1 = \diag(1,1,-1,-1)$. We can rescale the coordinates so that $m_1\colon x_0 - x_1 = x_2 - x_3 = 0$ and
\[
\tau_2 = \left[\begin{array}{cccc}   
&  & 1  & \\ 
 & & & 1  \\
1 & & & \\
 & 1 & &
\end{array}\right].
\]

If $l_0$ is given by $x_0 - \mu x_2 = x_1 - \nu x_3$, $\mu,\nu\in \CC$, then there exists $\lambda \in \CC$ so that $Q_{54}$ is given by an equation of the following form:
\begin{multline} \label{eq:Q54}
x_{0}^{4} + x_{2}^{4} + \lambda {\left(x_{1}^{4} + x_{3}^{4}\right)} - \frac{{\left(\mu^{4} + 1\right)}}{\mu^{2}} x_{0}^{2} x_{2}^{2} - \frac{\lambda {\left(\nu^{4} + 1\right)} }{\nu^{2}} x_{1}^{2} x_{3}^{2} \\ 
- {\left(\lambda + 1\right)}{\left(x_{0}^{2} x_{1}^{2} + x_{2}^{2} x_{3}^{2}\right)} + \frac{{\left(\mu^{2} \nu^{3} \lambda - \mu^{3} \nu^{2} - \mu \lambda + \nu\right)}}{{\left(\mu - \nu\right)} \mu \nu} {\left(x_{1}^{2} x_{2}^{2} + x_{0}^{2} x_{3}^{2}\right)}
\\
- \frac{{\left(\mu^{2} \nu^{4} \lambda - \mu^{4} \nu^{2} - \mu^{2} \lambda + \nu^{2}\right)} {\left(\mu + \nu\right)}}{{\left(\mu - \nu\right)} \mu^{2} \nu^{2}} x_{0} x_{1} x_{2} x_{3} = 0.
\end{multline}

We will explain below how to determine $\lambda,\mu$ and $\nu$. It turns out that
\begin{align*}
\lambda = & \frac{95}{432} \, \nu^{11} - \frac{49}{48} \, \nu^{10} + \frac{1019}{432} \, \nu^{9} - \frac{1447}{432} \, \nu^{8} - \frac{323}{216} \, \nu^{7} + \frac{2221}{216} \, \nu^{6} \\ & - \frac{1247}{216} \, \nu^{5} - \frac{485}{216} \, \nu^{4} + \frac{3997}{144} \, \nu^{3} - \frac{7345}{432} \, \nu^{2} + \frac{593}{144} \, \nu - \frac{23}{48},
\\
 \mu = & \frac{1}{96} \, \nu^{11} - \frac{13}{288} \, \nu^{10} + \frac{19}{288} \, \nu^{9} + \frac{1}{288} \, \nu^{8} - \frac{55}{144} \, \nu^{7} + \frac{113}{144} \, \nu^{6} \\ & + \frac{41}{144} \, \nu^{5} - \frac{205}{144} \, \nu^{4} + \frac{379}{288} \, \nu^{3} + \frac{59}{288} \, \nu^{2} - \frac{359}{96} \, \nu + \frac{1}{32},
\end{align*}
and the minimal polynomial of $\nu$ over $\QQ$ is
\begin{multline*}
\nu^{12} -6 \, \nu^{11} + 18 \, \nu^{10} - 34 \, \nu^{9} + 23 \, \nu^{8} + 44 \, \nu^{7} - 100 \, \nu^{6} \\ + 68 \, \nu^{5} + 127 \, \nu^{4} - 262 \, \nu^{3} + 242 \, \nu^{2} - 66 \, \nu + 9.
\end{multline*}

\begin{remark}
As a matter of fact, the surface $Q_{54}$ is defined over $\QQ(\lambda)$, which is a non-Galois field extension of degree~$6$ of $\QQ$. All lines are defined over its Galois closure $\QQ(\nu) = \QQ(\lambda,i)$.
\end{remark}

\subsubsection{An explicit isomorphism between \texorpdfstring{$Q_{54}$}{Q54} and \texorpdfstring{$X_{52}''$}{X52ii}}

Let $X_{52}''$ be the only surface up to projective equivalence containing configuration~$\bX_{52}''$. Since the transcendental lattices of $Q_{54}$ and $X_{52}''$ are in the same oriented isometry class, these surfaces are isomorphic to each other by a theorem of Shioda-Inose~\cite{shioda-inose77}. More precisely, they form an Oguiso pair (cf. Section~\ref{sec:oguiso-pairs}). In Section~\ref{subsec:X52ii} we explain how to find a defining equation of $X_{52}''$. Starting from this equation, we provide here a way to compute an explicit isomorphism between $Q_{54}$ and $X_{52}''$ following a method illustrated by Shimada and Shioda. We refer to their article~\cite{shimada-shioda} for further details on the algorithms used.

Let $(S,h)$ be the configuration $\bX_{52}''$. Let $\mathcal L$ be the set of the $52$ lines in $S$,
\[
\cL := \{ l \in S \mid l^2 = -2,\, l\cdot h = 1 \}.
\]
Compute the set of very ample polarizations of degree~$4$ which have intersection~$6$ with $h$,
\[
\cH := \{v \in S \mid v^2 = 4, v\cdot h = 6, \,\text{$v$ very ample}\}
\]
(cf. also~\cite[Lemma 6.8]{degtyarev16}). The set $\cH$ has $153$ elements.
Let $\cO$ be the set of vectors $v$ in $\cH$ such that
\begin{enumerate}
\item the configuration $(S,v)$ is isometric to $\bQ_{54}$;
\item there are six pairwise distinct lines $l_0,\ldots,l_5\in \cL$ such that \[ v = 3\,h - l_0 - \ldots - l_5. \]
\end{enumerate}
There are $36$ vectors in $\cH$ which satisfy the first condition. (The other $117$ vectors define a configuration isometric to $\bX_{52}''$.) The set $\cO$ has $6$ elements. For any vector $v\in \cO$ and sextuple $l_0,\ldots,l_5 \in \cL$ as above, it turns out that, up to relabeling, $l_0$ is of type $(4,4)$, $l_1,l_2$ are of type $(4,3)$, $l_3$ is of type $(3,5)$ and $l_4,l_5$ are of type $(0,12)$.

Fix a vector $v \in \cO$ and a sextuple $l_0,\ldots,l_5 \in \cL$ as above. Compute the explicit equations of the lines $l_i$ in the surface $X_{52}''$ (cf. Remark~\ref{rmk:X52ii-def-lines}). 

A defining equation of $Q_{54}$ can then be obtained in a similar way as in~\cite[Theorem 4.5]{shimada-shioda}.
Let $\Gamma_d$ be the space of homogeneous polynomials of degree $d$ in the variables $x_0,x_1,x_2,x_3$. Let $\Lambda \subset \Gamma_4$ be the 4-dimensional subspace of cubic polynomials that vanish along the lines $l_0,\ldots,l_5$. Since we know the equations of these lines, we can compute explicitly a basis $\varphi_0,\ldots,\varphi_3$ of $\Lambda$. Let $\bar\Gamma \subset \Gamma_{12}$ be the $290$-dimensional subspace of polynomials of degree~12 whose degree with respect to $x_0$ is $\leq 3$. 
Let $\sigma\colon \Gamma_4 \rightarrow \Gamma_{12}$ be the homomorphism given by the substitution $x_i \mapsto \varphi_i$. 
Let $\rho\colon \Gamma_{12} \rightarrow \bar\Gamma_{12}$ be the homomorphism given by the remainder of the division by the defining polynomial \eqref{eq:X52ii} of $X_{52}''$. Then, the kernel of $\rho\circ\sigma$ has dimension~1 and is generated by a defining equation of~$Q_{54}$.

It is then a matter of changing coordinates in order to find an equation as in~\eqref{eq:Q54}.

\begin{remark} Let $\zeta$ be a primitive 3rd root of unity. Let $r$ be the algebraic number defined as in Remark~\ref{rmk:X52ii-def-lines}.
The isomorphism that we found is defined over a degree~24 Galois extension of $\QQ$ generated by the element $x = ir$.  Its minimal polynomial is
\begin{multline*}
x^{24} + 38 \, x^{22} + 1045 \, x^{20} + 16306 \, x^{18} + 180538 \, x^{16} - 258514 \, x^{14} + 166541 \, x^{12} \\ - 258514 \, x^{10} + 180538 \, x^{8} + 16306 \, x^{6} + 1045 \, x^{4} + 38 \, x^{2} + 1.
\end{multline*}
Note that $\QQ(x) = \QQ(\nu,\zeta) = \QQ(r,i)$.
\end{remark}


\section{Some configurations with at most 52 lines}
\label{sec:<=52-lines}

In this section we give explicit equations of the surfaces containing configurations $\bX_{52}''$, $\bX_{52}'''$, $\bX_{52}^\mathrm v$, $\bY_{52}'$, $\bQ_{52}'''$ and $\bX_{51}$. The same conventions as in §\ref{subsec:conventions} apply.

\subsection{Configuration \texorpdfstring{$\bX_{52}''$}{X52ii}}
\label{subsec:X52ii}

Let $X_{52}''$ be a quartic surface containing configuration $\bX_{52}''$. Let $l_0$ be the only line of type~$(6,0)$ contained in $X_{52}''$. There is a symplectic automorphism~$\sigma$ of order~$3$ which preserves $l_0$ and six of its reducible fibers. By Proposition~\ref{prop:special-aut}, $l_0$ is a special line; hence, $X_{52}''$ is given by an equation as in \eqref{eq:familyZ} and $\sigma = \diag(1,1,\zeta,\zeta^2)$, with $\zeta$ a primitive 3rd root of unity. 

Let $m\colon x_2 = x_3 = 0$ be the line in $\PP^3$ fixed pointwise by $\sigma$. Let $\tau$ be a symplectic automorphism of $X_{52}''$ of order~$2$. Then $\tau$ fixes four lines $l_0,\ldots,l_3$ and $\sigma^2\tau = \tau\sigma$, which implies that $\tau(\Fix(\sigma)) = \Fix(\sigma)$. The line $l_0$ has two ramified fibers in $\Pi'\colon x_0 = 0$ and $\Pi''\colon x_1 = 0$, so the planes $\Pi',\Pi''$ are permuted by $\tau$. If $\Pi',\Pi''$ were fixed by $\tau$, then $\tau$ would commute with $\sigma$. It follows that $\tau(\Pi') = \Pi''$, so $\tau$ has the following form (after rescaling one variable):
\begin{equation} \label{eq:X52ii-tau-a}
\tau = \left[\begin{array}{cccc}   
  	& p &   & \\ 
 1  &  	& 	&   \\
 	&  	& 	& p \\
 	&  	& 1	&  \\ 
\end{array}\right],
\end{equation}
for some $p \in \CC$. For $i=1,2,3$, we can suppose that $l_i$ is the intersection of $\Pi_i\colon e_i x_0 + c_i x_1 + d_i x_2 = 0$ and $\tau(\Pi_i)$. After rescaling, we can suppose that $c_1 = d_1 = e_1 = 1$; we let $c := c_2$, $d := d_2$, $e := e_2$. 

Imposing that $l_1$ and $l_2$ are contained in $X_{52}''$ we can express all coefficients in terms of $p,c,d,e$; moreover, one of the following equations must hold:
\begin{align}
p &= \frac{{\left(c - d\right)}^{2}}{{\left(d - e\right)}^{2}}, \label{eq:X52ii-comp3} \\
p &= \frac{c^{2} + c d + d^{2}}{d^{2} + d e + e^{2}}, \label{eq:X52ii-comp2} \\
p &= -\frac{{\left(c + d\right)}^{2} {\left(c - d\right)}}{{\left(d^{2} - c e\right)} {\left(d + e\right)}}. \label{eq:X52ii-comp1}
\end{align}
Condition \eqref{eq:X52ii-comp3} implies that $l_1$ and $l_2$ intersect each other. Condition \eqref{eq:X52ii-comp2} implies that there exists a line intersecting $l_0$, $l_1$, $l_2$ and $l_3$. Both are contradictions, so condition~\eqref{eq:X52ii-comp1} must hold. 

We parametrize the pencil of planes containing $l_1$ by
\[
t \mapsto \{x_{2} = -x_{0} - x_{1} - t {\left(x_{1} + x_{3} + x_{0}/p\right)}\}.
\]
The discriminant $\Delta$ of $l_1$ has the following form:
\[
\Delta = P Q^2 R^2,
\]
where $P$, $Q$, $R$ are polynomials in $t$ of degree $12$, $4$, $2$, respectively. The line $l_1$ and $l_2$ are of type $(4,4)$, so $R$ divides $P$ and the following condition holds:
\begin{multline} \label{eq:X52ii-genus1}
c^{3} d^{2} + c^{2} d^{3} + c d^{4} - c^{4} e + 2 \, c^{3} d e + 4 \, c^{2} d^{2} e + 2 \, c d^{3} e \\ - d^{4} e + c^{3} e^{2} + c^{2} d e^{2} + c d^{2} e^{2} = 0.
\end{multline}

Under condition \eqref{eq:X52ii-genus1}, the polynomial $Q$ splits into two degree~2 polynomials $Q = Q_1Q_2$, with
\begin{multline*}
Q_2 = \left( c^{2} d^{3} + c d^{4} + d^{5} - c^{3} d e + d^{4} e - c^{3} e^{2} - c^{2} d e^{2} - c d^{2} e^{2}\right) t^{2} \\ 
+ \left(- c^{3} d^{2} + c^{2} d^{3} + 4 \, c d^{4} + 2 \, d^{5} + c^{4} e + 2 \, c^{3} d e + 2 \, c^{2} d^{2} e + c d^{3} e \right) t \\
- c^{5} - 2 \, c^{4} d - c^{3} d^{2} + c^{2} d^{3} + 2 \, c d^{4} + d^{5};
\end{multline*}
moreover, $P = W Q_2 R$ for some polynomial $W$ of degree~8, which we can explicitly compute. The polynomial $W$ has two double roots which account for the remaining $2$-fibers of $l_1$. After normalizing $d$ to $1$, we compute the resultant of $W$, obtaining another condition on $c,e$. Together with \eqref{eq:X52ii-genus1}, we get that
\[
e = -\frac{192}{247} \, c^{5} - \frac{60}{247} \, c^{4} - \frac{3624}{247} \, c^{3} - \frac{621}{19} \, c^{2} - \frac{5336}{247} \, c - \frac{708}{247},
\]
and $c$ satisfies
\[
c^6 + 19\,c^4 + 36\,c^3 + 19\,c^2 + 1 = 0.
\]

\begin{remark}
The last computations can be simplified by noting that all factors $F$ of $\Delta$ satisfy a condition of symmetry due to $\tau$:
\[
 k^{n/2} F(t) = t^n F(k/t),
\]
where $k = (c+d)^2(d-c)/((d^2-ce)(d+e))$ and $n = \deg F$.
\end{remark}

\begin{remark}
The field $\QQ(c)$ is a Galois extension of $\QQ$ of degree~6. The surface $X_{52}''$, though, can be defined on the smaller non-Galois extension $\QQ(p)$, where $p$ is the parameter appearing in \eqref{eq:X52ii-tau-a} and is equal to
\[
p = \frac{3}{26} c^{5} + \frac{3}{13} c^{4} + \frac{12}{13} c^{3} + \frac{15}{13} c^{2} + \frac{15}{26} c + \frac{4}{13}.
\]
The minimal polynomial of $p$ over $\QQ$ is 
\[ p^{3} - 201 \, p^{2} + 111 \, p - 19. \]

The surface $X_{52}''$ is then defined by
\begin{multline} \label{eq:X52ii}
x_{0}^{4} + \left(-36 p^{2} + 6696 p + 2052\right) x_{0}^{3} x_{1} + 4968 p x_{0}^{2} x_{1}^{2} \\ + \left(-540 p^{2} + 6048 p - 684\right) x_{0} x_{1}^{3} + 3312 p^{2} x_{1}^{4} + \left(-19 p^{2} - 100 p + 209\right) x_{1} x_{2}^{3} \\ + \left(11 p^{2} - 5542 p + 1121\right) x_{0}^{2} x_{2} x_{3} + \left(-116 p^{2} + 1612 p - 380\right) x_{0} x_{1} x_{2} x_{3} \\ + \left(-3331 p^{2} - 100 p + 209\right) x_{1}^{2} x_{2} x_{3} + \left(-3919 p^{2} + 2318 p - 361\right) x_{0} x_{3}^{3} = 0.
\end{multline}

\begin{remark} \label{rmk:X52ii-def-lines}
All lines on $X_{52}''$ are defined over $K = \QQ(c,\zeta)$, where $\zeta$ is a primitive 3rd root of unity. A primitive element of $K$ over $\QQ$ is $r = c\zeta$, whose minimal polynomial is 
\[
r^{12} - 19 \, r^{10} + 72 \, r^{9} + 342 \, r^{8} - 684 \, r^{7} + 937 \, r^{6} - 684 \, r^{5} + 342 \, r^{4} + 72 \, r^{3} - 19 \, r^{2} + 1.
\]

\end{remark}

\end{remark}

\begin{remark}
The quintic curve in $\PP^2$ defined by condition~\eqref{eq:X52ii-genus1} has one singularity of type $\bD_4$ and two singularities of type $\bA_1$. Hence, it has geometric genus~$1$; in particular, it is not rational.
\end{remark}

\subsection{Configuration \texorpdfstring{$\bX_{52}'''$}{X52iii}}
\label{subsec:X52iii}

Let $X_{52}'''$ be a quartic surface containing configuration $\bX_{52}'''$. 
Then $X_{52}'''$ contains a pair of twin lines, say $l_0$ and $l_1$, of type $(0,10)$. There are four symplectic symmetries of order~$5$, which fix both $l_0$ and $l_1$, and no other lines. Choose one of them and call it $\tau$. By~§\ref{subsec:sympl-autom}, $\tau$ has exactly $4$ fixed point on $X_{52}'''$. Necessarily, two of them lie on $l_0$ and the other two on $l_1$. 

Up to coordinate change, we can suppose that $l_0$ and $l_1$ are given by $x_0 = x_1 = 0$ and $x_2 = x_3 = 0$, respectively, and that the fixed points of $\tau$ are the coordinate points. 
It follows that $\tau = \diag(1,\xi,\xi^i,\xi^j)$, where $\xi$ is a primitive 5th root of unity, $0\leq i,j\leq 4$. 
In fact, the first and second entries cannot be equal, because $\tau$ does not fix $l_0$ pointwise. Analogously for $l_1$, we have $i \neq j$; hence, up to exchanging $x_2$ with $x_3$, we can suppose $i < j$. 

The conditions $i = 0$ or $i = 1$ lead to contradictions: $X_{52}'''$ would contain more lines fixed by $\tau$ than just $l_0$ and $l_1$. As $X_{52}'''$ is smooth, we find $i = 2$ and $j = 4$. 

Imposing that $\tau$ is a symplectic symmetry and normalizing the remaining coefficients, $X_{52}'''$ turns out to be a Delsarte surface:
\[
x_{0}^{3} x_{2} + x_{1} x_{2}^{3} + x_{1}^{3} x_{3} + x_{0} x_{3}^{3} = 0.
\]

All lines intersecting both $l_0$ and $l_1$ (for instance, the line given by $x_0 - x_1 = x_2 + x_3 = 0$) are of type $(4,6)$.

\subsection{Configuration \texorpdfstring{$\bX_{52}^\mathrm v$}{X52v}}
\label{subsec:X52v}

Let $X_{52}^\mathrm v$ be a quartic surface whose Fano configuration is isometric to $\bX_{52}^\mathrm v$. Then $X_{52}^\mathrm v$ contains $12$ lines of type $(2,8)$, but only four of them, say $l_0,\ldots,l_3$, form a quadrangle in which the opposite lines form two pairs of twin lines. We choose coordinates so that $l_i$ is the line $x_i = x_{i+1} = 0$.

There is a symplectic symmetry $\sigma$ of order~$4$ such that $\sigma(l_0) = l_2$, $\sigma(l_1) = l_3$ and $\sigma^2$ fixes $l_0,\ldots,l_3$. Since $\sigma^2$ has exactly 8 fixed points, of which 4 are the vertices of the quadrangle, we have $\sigma^2 = \diag(1,-1,1,-1)$. It follows that $\sigma$ is given by
\[
\left[\begin{array}{cccc}   
 	&  	& 1 &  	\\ 
 	&  	&  	& a \\
-ab &  	&  	& 	\\
    & b & 	& 
\end{array}\right]
\]
for some $a,b \in \CC$. The residual conics in the coordinate planes are irreducible. Hence, there is a plane containing $l_0$ different from $x_0 = 0$ and $x_1 = 0$ where the residual cubic splits into three lines. We let $m_1$ be the line intersecting $l_2$, which must be of type $(5,3)$. By Corollary~\ref{cor:p-fiber-twin-line}, the point of intersection of the other two lines lies on $l_0$. We introduce two parameters $p,q \in \CC$, so that $m_1$ is given by $x_0 - px_1 = x_2 - qx_3 = 0$. After imposing all conditions, we normalize all remaining coefficients except $a$ by rescaling variables. We find that $X_{52}^\mathrm v$ is given by a polynomial of the form
\begin{multline*}
a^{4} x_{1} x_{3}^{3} - a^{3} x_{1}^{3} x_{3} - a x_{0}^{3} x_{2} - {\left(a^{3} - 2 \, a\right)} x_{0} x_{1}^{2} x_{2} + {\left(2 \, a^{3} - a\right)} x_{0}^{2} x_{1} x_{3} \\ - {\left(2 \, a^{2} - 1\right)} x_{1} x_{2}^{2} x_{3} - {\left(a^{4} - 2 \, a^{2}\right)} x_{0} x_{2} x_{3}^{2} - x_{0} x_{2}^{3} = 0.
\end{multline*}

In order to determine $a$, we inspect the discriminant $\Delta$ of the fibration induced by $m_1$. We parametrize the planes containing $m_1$ by
\[
t \mapsto \{ x_0 = px_1 + t(x_2 - qx_3)\}.
\]
It turns out that the discriminant is of the form
\[ \Delta = t^4 P Q^2 R^3, \]
where $P$, $Q$ and $R$ are polynomials in $t$ of degree $4$, $4$ and $2$, respectively. Since $m_1$ has type $(5,3)$, $P$ and $Q$ must have a common root. We compute the determinant of the Sylvester matrix associated to $P$ and $Q$, finding that $a$ must be a root of one of the following polynomials
\begin{align*}
& a^{4} - 3 \, a^{3} + 2 \, a^{2} + 3 \, a + 1, \\
& a^{4} + 3 \, a^{3} + 2 \, a^{2} - 3 \, a + 1, \text{ or} \\
& a^{12} - 4 \, a^{10} + 2 \, a^{8} + 5 \, a^{6} + 2 \, a^{4} - 4 \, a^{2} + 1.
\end{align*}
All solutions lead to projective equivalent surfaces, as $\Cl(\bX_{52}^\mathrm v) = 1$, see §\ref{subsec:proj-equiv-classes}.

\subsection{Configuration \texorpdfstring{$\bY_{52}'$}{Y52i}} \label{subsec:Y52i}
Let $Y_{52}'$ be a quartic surface containing configuration $\bY_{52}'$. Then $Y_{52}'$ contains two intersecting lines $l_0$ and $l_1$ of type $(4,6)$. We let $l_0'$ and $l_1'$ be their respective twin lines. As these lines form a quadrangle, we can choose coordinate so that $l_0\colon x_0 = x_1 = 0$, $l_0'\colon x_2 = x_3 = 0$, $l_1\colon x_0 = x_3 = 0$, $l_1'\colon x_1 = x_2 = 0$. The residual conics in the coordinate planes are irreducible. Moreover, there exists a unique symplectic symmetry $\sigma$ of order $2$ which fixes these four lines. It follows that the two lines fixed in $\PP^3$ by $\sigma$ must be the only two lines not contained in $Y_{52}'$ joining two vertices of the quadrangle, namely $x_0 = x_2 = 0$ and $x_1 = x_3 = 0$, so $\sigma = \diag(1,-1,1,-1)$. 

Let $\tau$ be one of the two symplectic symmetries of order $4$. We have $\tau(l_0) = l_1$ and $\tau(l_0') = l_1'$ and $\tau^2 = \sigma$. After rescaling variables, $\tau$ is given by
\[
\tau = \left[\begin{array}{cccc}   
1 &  &   &  \\ 
 &   &  & 1  \\
 &  & r & \\
 & -1 & &
\end{array}\right],
\]
with $r^2 = 1$. If $r = 1$, then the conic in $x_0 = 0$ is reducible, which is not the case, so $r = -1$. Hence, $Y_{52}'$ is given by an equation of the form
\[
a x_{0}^{3} x_{2} + b x_{0} x_{2}^{3} + c x_{0}^{2} x_{1} x_{3} + d x_{1} x_{2}^{2} x_{3} + x_{0} x_{1}^{2} x_{2} + x_{1}^{3} x_{3} + x_{0} x_{2} x_{3}^{2} + x_{1} x_{3}^{3} = 0,
\]
for some $a,b,c,d \in \CC$. This is a 3-dimensional family, as one more coefficient can be set to $1$ by rescaling the other coefficients and the variables.

We introduce two new parameters $p,q \in \CC$ so that one of the lines $m$ of type $(3,5)$ intersecting both $l_0$ and $l_0'$ and contained in a $3$-fiber of $l_0$ is given by $m\colon x_0 - px_1 = x_2 - qx_3 = 0$.  By Corollary~\ref{cor:p-fiber-twin-line}, the other two lines in the plane containing $l_0$ and $m$ intersect in a point lying on $l_0$. The line $\tau(m)$ intersects $m$, and their residual conic is also reducible.

Imposing all these conditions and normalizing $q$, we can express all coefficients in terms of $p$.  We find that
\begin{align*}
 a & = -\frac{{\left(2 \, p + 1\right)} {\left(p + 1\right)}^{2}}{2 \, {\left(p^{2} - 2 \, p - 1\right)} p^{3}}, \\
 b & = -\frac{p^{2} + 2 \, p + 1}{4 \, p}, \\
 c & = \frac{{\left(7 \, p + 3\right)} {\left(p + 1\right)}}{2 \, {\left(p^{2} - 2 \, p - 1\right)} p^{2}}, \\
 d & = \frac{1}{4} \, {\left(p + 1\right)} {\left(p - 3\right)}.
\end{align*}

Finally, we require that $l_0$ has type $(4,6)$ by looking at the discriminant $\Delta$ of its induced fibration. We parametrize the planes containing $l_0$ by $t \mapsto \{x_0 = tx_1\}$. The discriminant has the following form:
\[
\Delta = t^2(t-p)^3(t+p)^3 P^2 Q^2 R^2 S,
\]
where $P,Q,R,S$ are polynomials in $t$ of degree $2$, and $P(-t) = Q(t)$. Therefore, $S$ must have a common root with either $P$, $Q$ or $R$. We find a finite list of values for $p$. Looking also at the determinant of $m$, it turns out that $p$ is a root of 
\[
p^{3} - \frac{11}{9} \, p^{2} - \frac{7}{3} \, p - 1.
\]

Since $\Cl(\bY_{52}') = 3$, each root corresponds to a different projective equivalence class. In particular, the real root corresponds to the surface with transcendental lattice $[2,0,38]$.

\subsection{Configuration \texorpdfstring{$\bQ_{52}'''$}{Q52iii}} \label{subsec:Q52iii}

Let $Q_{52}'''$ be a quartic surface containing configuration $\bQ_{52}'''$. Then $Q_{52}'''$ contains exactly 4 lines $l_0,\ldots,l_3$ of type $(5,0)$. Since these lines intersect each other pairwise, they are coplanar. Moreover, there exists a symplectic automorphism $\sigma$ which fixes each~$l_i$. It follows that $l_0,\ldots,l_3$ form a star, otherwise $\sigma$ would fix three points on at least one $l_i$, so it would fix the whole line, but this cannot happen, as $\sigma$ is symplectic. By Lemma~\ref{lem:autP3-ord2}, $\sigma$ fixes two lines $m',\,m''$ in $\PP^3$ pointwise. If $P$ is the center of the star and $Q_i$ is the other point on $l_i$ fixed by $\sigma$ for $i=0,\ldots,3$, then necessarily all $Q_i$ lie on one of the two lines, say $m'$, and the other line $m''$ contains $P$. All symmetries of $Q_{52}'''$ fix $P$.

There is a non-symplectic symmetry $\tau$ of order 3 which fixes each $l_i$ and commutes with $\sigma$; hence, $\tau$ fixes each $Q_i$. This means that, as an automorphism of $\PP^3$, $\tau$ fixes $m'$ pointwise. Its invariant lattice has rank 4, so by results of Artebani, Sarti and Taki \cite{ArtebaniSartiTaki}, $\Fix(\tau,X)$ consists of one smooth curve $C$ and exactly one point. By Lemma~\ref{lem:autPn-ord3}, the curve $C$ must be the intersection of $Q_{52}'''$ with a plane $\Pi$ in $\Fix(\tau,\PP^3)$; moreover, $\Pi$ necessarily contains $m'$, but not $P$, since $C$ is smooth. Let $R$ be the point of intersection of $\Pi$ and $m''$. If $R\in X$, then $\tau$, being of order 3, would fix at least three points on $m''$, so it would fix the whole line $m''$ pointwise. This is impossible since $\Fix(\tau,\PP^3) = \Pi \cup \{P\}$, so $R\notin X$ and $m''\cap X$ consists of four distinct points.

There exists also a symplectic symmetry $\varphi$ such that $\varphi^2 = \sigma$. Since $\varphi$ commutes with $\sigma$ and $\tau$, $\varphi(R) = R$, so $\varphi$ fixes $m''$ pointwise. Up to relabeling, $\varphi(l_0) = l_3$ and $\varphi(l_1) = l_2$. We choose coordinates in such a way that $P = (0, 0, 1, 0)$, $Q_0 = (0, 0, 0, 1)$, $Q_3 = (0, 1, 0, 0)$, and $R = (1, 0, 0, 0)$.
%
%
%

After rescaling, the automorphisms $\tau$ and $\varphi$ are given by
\[
\tau = \left[\begin{array}{cccc}   
1&  &   &  \\ 
 & 1 &  &   \\
 &  & \zeta & \\
 & & & 1 
\end{array}\right] \quad \text{and} \quad
\varphi = \left[\begin{array}{cccc}   
1&  &   &  \\ 
 &  &  & -1  \\
 &  & 1 & \\
 & 1 & & 
\end{array}\right],
\]
where $\zeta$ is a primitive 3rd root of unity. If $\omega$ is a non-zero 2-form on $Q_{52}'''$, then either $\tau^*\omega = \zeta\omega$ or $\tau^*\omega = \zeta^2\omega$, but the second condition leads to $P$ being a singular point. 

It follows that $Q_{52}'''$ is given by an equation of the form
\[
c x_{0}^{2} x_{1}^{2} + c x_{0}^{2} x_{3}^{2} + d x_{1}^{2} x_{3}^{2} + x_{0}^{4} + x_{0} x_{2}^{3} - x_{1}^{3} x_{3} + x_{1} x_{3}^{3} = 0,
\]
for some $c,d\in \CC$. Parametrizing the planes containing $l_0$ by $t \mapsto \{x_0 = tx_1\}$, the discriminant $\Delta$ of the fibration induced by $l_0$ has the following form:
\[
\Delta = t^4 P(t^2)^2,
\]
where $P$ is a polynomial of degree $5$. Looking at the resultant of $P$ and excluding the conditions leading to surface singularities, we find that
\[
27 \, c^{4} d^{2} - 54 \, c^{2} d^{3} + 100 \, c^{4} + 27 \, d^{4} - 198 \, c^{2} d + 162 \, d^{2} + 243 = 0.
\]
This polynomial splits over $\QQ(\zeta)$. We see that there exists $e \in \CC$ such that
\[
c = \frac{e^2-2\,\zeta-1}{3\,e}, \quad \text{and} \quad d = \frac 19 (e^2-20 \, \zeta - 10).
\]

Moreover,
\[
\Delta = t^4 (t^2 - e/3)^4 Q(t^2)^2,
\]
where $Q$ is a polynomial of degree 3. Looking now at the resultant of $Q$, it turns out that if $e$ is a root of
\[
e^{4} - 20 \, {\left(2 \, \zeta + 1\right)} e^2 + \frac{15}{4},
\]
then $l_0$ is of type $(5,0)$, so the Fano configuration of $Q_{52}'''$ is indeed $\bQ_{52}'''$. As $\Cl(\bQ_{52}''') = 1$, all other quartic surfaces with this Fano configuration are projectively equivalent to the one we found.

\subsection{Configuration \texorpdfstring{$\bX_{51}$}{X51}} \label{subsec:X51}

Let $X_{51}$ be a surface containing configuration~$\bX_{51}$. Then $X_{51}$ contains a line $l_0$ of type $(6,2)$. By Proposition~\ref{prop:(6,q)-is-special}, $X_{51}$ is given by an equation as in~\eqref{eq:familyZ}. In particular, the two lines $l_1$ and $l_2$ in the $1$-fibers are given by $x_0 = x_2 = 0$ and $x_1 = x_3 = 0$.
Note that the residual conics in the $1$-fibers intersect $l_1$ and $l_2$ at the coordinate points and are tangent to~$l_0$. 

There are three symplectic symmetries of order $2$ which exchange $l_1$ and $l_2$. Choose one of them and call it $\sigma$; necessarily, $\sigma$ has the following form:
\begin{equation*} 
\sigma = \left[\begin{array}{cccc}   
   &   1 &   & \\
ab &&& \\
 &&&   a \\
 &&   b   & \\
\end{array}\right].
\end{equation*}
Moreover, $\sigma$ fixes one line $m$ intersecting both $l_1$ and $l_2$. Up to rescaling variables, $m$ is given by $x_3 - x_1 = x_2 - ax_0 = 0$. By further inspection of $\bX_{51}$, we see that the residual conic in the plane containing $m$ and $l_1$ is reducible. After normalizing all coefficients except $a$, we are left with the following $1$-dimensional family:
\begin{multline*}
3 \, a^{3} x_{0}^{2} x_{1}^{2} - 3 \, a^{3} x_{0}^{2} x_{2} x_{3} - 3 \, a^{2} x_{0} x_{1} x_{2} x_{3} - 3 \, a^{2} x_{1}^{2} x_{2} x_{3} \\ - {\left(a^{4} - a^{3}\right)} x_{0}^{3} x_{1} - {\left(a^{3} - a^{2}\right)} x_{0} x_{1}^{3} + {\left(4 \, a - 1\right)} x_{1} x_{2}^{3} + {\left(4 \, a^{3} - a^{2}\right)} x_{0} x_{3}^{3} = 0.
\end{multline*}

In order to find the last condition for $a$, we investigate the discriminant~$\Delta$ of the fibration induced by $l_1$. We parametrize the planes containing $l_1$ by $t\mapsto \{x_0 = tx_2\}$, so that $\Delta$ has the following form:
\[
\Delta = t^3 (a^3t^3 -1)^3 P,
\]
where $P$ is a polynomial in $s = t^3$ of degree~$4$. As the line $l_1$ is of type $(3,4)$, the resultant of $P$ must vanish. Knowing also that $m$ is of type $(2,7)$, it follows that $a$ satisfies the following equation:
\begin{equation*}
a^{3} - \frac{11}{3} \, a^{2} + \frac{10}{9} \, a - \frac{1}{9} = 0.
\end{equation*}
Since $\Cl(\bX_{51}) = 3$, each root corresponds to a different projective equivalence class. In particular, the real root corresponds to the surface with transcendental lattice $[6,3,16]$.

\subsection{Configuration \texorpdfstring{$\bZ_{52}$}{Z52}} \label{subsec:Z52}

The general member of the following rational family contains $52$ lines forming configuration~$\bZ_{52}$ and has Picard number 19:
\begin{multline*} 
    t^{2} x_{1} x_{2} {\left(t x_{0} + t x_{3} - 2 \, x_{1} + 2 \, x_{2}\right)} {\left(t x_{0} - a x_{3} - 2 \, x_{1} - 2 \, x_{2}\right)}  \\
    = -4 \, x_{0} x_{3} {\left(t x_{0} + t x_{3} - 6 \, x_{1} + 6 \, x_{2}\right)} {\left(t x_{0} - t x_{3} - 6 \, x_{1} - 6 \, x_{2}\right)}.
\end{multline*}

This family was found by taking advantage of the fact that a surface containing $\bZ_{52}$ has four special lines of type~$(6,0)$ and six pairs of twin lines of type~$(2,8)$.

Generically, the lines $x_0 = x_1 = 0$ and $x_2 = x_3 = 0$ are twin lines, while the lines $x_0 = x_2 = 0$ and $x_1 = x_3 = 0$ are special lines. All surfaces of the family have the symmetry
\[
\left[\begin{array}{cccc}   
    &     & &1 \\
  	& & -1 & \\
 	& 1 && \\
 -1 & &  & \\
\end{array}\right].
\]

We obtain models containing configurations $\bX_{64}$ and $\bX'_{60}$ when the minimal polynomial of $t$ is $t^{4} + 144 $ or $t^{4} - 12 \, t^{2} + 144$, respectively.

\section{An explicit Oguiso pair}
\label{sec:oguiso-pairs}

In this section we answer a question posed by Oguiso~\cite{oguiso17}.

\subsection{Oguiso pairs}
\label{subsec:oguiso-pairs}
Let $\pi_1,\pi_2\colon \PP^3\times \PP^3\rightarrow \PP^3$ be the first and second projection. A pair $(X_1,X_2)$ of smooth quartic surfaces (not necessarily distinct) is an \emph{Oguiso pair} if there exists a smooth complete intersection $S$ of four hypersurfaces $Q_0,\ldots,Q_3$ of bi-degree $(1,1)$ in $\PP^3\times \PP^3$ such that the restriction to $S$ of $\pi_i$ is an isomorphism onto $X_i$, for $i=1,2$. In particular, $X_1$ and $X_2$ are isomorphic as abstract K3 surfaces, and the isomorphism is given by 
\[
(\pi_2|_S) \circ (\pi_1|_S)^{-1} \colon X_1 \xrightarrow{\sim} X_2.
\]

Conversely, let $X$ be a K3 surface and suppose $h_1, h_2$ are very ample divisors which induce embeddings of $X$ into $\PP^3$ whose images are $X_1$ and $X_2$, respectively.
\begin{theorem}[Oguiso \cite{oguiso17}] The smooth quartic surfaces $(X_1,X_2)$ form an Oguiso pair if and only if $h_1\cdot h_2 = 6$.
\end{theorem}

Under this assumption, both $X_1$ and $X_2$ are determinantal quartic surfaces (also called Cayley quartic surfaces). A determinantal description can be given in the following way. Let $x_0,\ldots,x_3$ and $y_0,\ldots,y_3$ be the coordinates in the first and second factor of $\PP^3\times\PP^3$, respectively. Write
\[
Q_{k} = \sum_{i,j,k=0}^3 a_{ijk} x_i y_j.
\]
Consider the matrix $M_1,M_2$ whose $(i,j)$-component are 
\begin{align*}
(M_1)_{ij} & = a_{0ij} x_0 + a_{1ij} x_1 + a_{2ij} x_2 a_{3ij}x_3, \\
(M_2)_{ij} & = a_{i0j} y_0 + a_{i1j} y_1 + a_{i2j} y_2 a_{i3j}y_3.
\end{align*}
Then the equations $\det M_1 = 0$ and $\det M_2 = 0$ define $X_1$ and $X_2$, respectively.

\subsection{Models of the Fermat quartic surface} Let $X_{48}$ be the Fermat quartic surface (which is the only surface up to projective equivalence containing configuration $\bX_{48}$) and let $X_{56}$ be one of the two surfaces containing configuration $\bX_{56}$ (which are complex conjugate to each other). Shioda first noticed that $X_{48}$ and $X_{56}$ are isomorphic to each other as abstract K3 surfaces. An explicit equation of $X_{56}$ and an explicit isomorphism between $X_{48}$ and $X_{56}$ were found by Shimada and Shioda~\cite{shimada-shioda}. The two surfaces are not projectively equivalent to each other, but they form an Oguiso pair. 

According to Degtyarev~\cite{degtyarev16}, there are no other smooth quartic models of the Fermat quartic and---curiously enough---$(X_{56},X_{56})$ is also an Oguiso pair, but $(X_{48},X_{48})$ is not (cf. \cite[§6.5]{degtyarev16}).

\subsection{A determinantal presentation} Oguiso asked in his paper~\cite{oguiso17} for explicit equations defining the complete intersection $S$ in $\PP^3\times\PP^3$ projecting onto $X_{48}$ and $X_{56}$. Here we provide such equations. In what follows we let $\zeta$ be a primitive 8th root of unity.

The explicit isomorphism $f\colon X_{48} \xrightarrow{\sim} X_{56}$ can be found in \cite[Table 4.1]{shimada-shioda}. The surface $S$ is the graph of $f$.

Generically, a hypersurface of bi-degree $(1,1)$ in $\PP^3\times \PP^3$ is defined by 16 coefficients. Choosing $12$ closed points $x_{1},\ldots,x_{12}$ on $X_{48}$ in a suitable way, one can find $12$ linearly independent conditions on these coefficients by imposing that $(x_i,f(x_i))$ belongs to $S$ for each $i=1,\ldots,12$.
The points chosen are 
\begin{align*}
(0,0,1,\zeta),\,(0,0,1,\zeta^5),\,(0,1,0,\zeta^5),\,(0,1,0,\zeta^7),\,(1,0,0,\zeta),\,(1,0,0,\zeta^5),\\
(1,0,0,\zeta^7),\,(1,0,\zeta^5,0),\,(1,\zeta,0,0),\,(0,1,\zeta,0),\,(0,1,\zeta^3,0),\,(\zeta,\zeta^2,\zeta,1),
\end{align*}
but of course this choice is quite arbitrary.

One ends up with a vector space of dimension $4$ of polynomials of bi-degree~$(1,1)$. The following four polynomials $Q_0,\ldots,Q_3$ form a basis of that vector space.

\begin{align*}
Q_0 ={}& \zeta^{3} x_{0} y_{2} + \zeta x_{3} y_{0} + x_{1} y_{0} - x_{2} y_{2}, \\
Q_1 ={}& \zeta^{3} x_{3} y_{1}  - \zeta^{3} x_{0} y_{3} + \zeta^{2} x_{2} y_{3} + x_{1} y_{1},
 \\
\begin{split}
Q_2 ={}& -{\left(\zeta^{2} - \zeta - 1\right)} x_{2} y_{0} - {\left(\zeta^{2} + 2 \, \zeta - 1\right)} x_{0} y_{1} - {\left(\zeta^{2} - 2\right)} x_{1} y_{1} \\
	& + {\left(\zeta^{2} - \zeta - 1\right)} x_{2} y_{1} + {\left(\zeta^{2} - 1\right)} x_{0} y_{2} - {\left(\zeta^{3} + \zeta - 1\right)} x_{1} y_{2} \\
    & - {\left(\zeta^{3} + \zeta + 2\right)} x_{3} y_{2} - {\left(\zeta^{3} + \zeta\right)} x_{0} y_{3} - {\left(\zeta^{3} + \zeta - 1\right)} x_{1} y_{3} \\
    & + {\left(\zeta^{2} - \zeta + 1\right)} x_{2} y_{3} +  \zeta x_{2} y_{2} + x_{1} y_{0},
\end{split} \\
\begin{split}
Q_3 ={}& 6
 \,\zeta^{3} x_{3} y_{3} - {\left(4 \, \zeta^{3} + 3 \, \zeta^{2} - 2 \, \zeta + 1\right)} x_{1} y_{0} - 3
 \, {\left(\zeta^{3} + \zeta\right)} x_{2} y_{0} \\
	& + {\left(2 \, \zeta^{3} + \zeta^{2} + 1\right)} x_{1} y_{1} - 3 \, {\left(\zeta^{3} + \zeta\right)} x_{2} y_{1} + 2 \, {\left(\zeta^{2} + \zeta - 1\right)} x_{0} y_{2} \\
	& + 3 \, {\left(\zeta^{2} + 1\right)} x_{1} y_{2} + {\left(3 \, \zeta^{3} + 2 \, \zeta^{2} - \zeta - 2\right)} x_{2} y_{2} - 2 \, {\left(\zeta^{3} + \zeta + 1\right)} x_{0} y_{3} \\
	& - 3 \, {\left(\zeta^{2} + 1\right)} x_{1} y_{3} + {\left(\zeta^{3} + \zeta - 2\right)} x_{2} y_{3} + 6\, x_{0} y_{0}
.\end{split} \\
\end{align*}

Indeed, the four hypersurfaces given by $Q_i = 0$ define a smooth complete intersection $S$ in $\PP^3\times \PP^3$ such that the restriction of $\pi_1$ and $\pi_2$ to $S$ is an isomorphism onto $X_{48}$ and $X_{56}$, respectively. This can be checked explicitly by computing the determinantal presentation explained in §\ref{subsec:oguiso-pairs}. The equation of $X_{56}$ that one obtains is the one provided in \cite[Theorem 1.3]{shimada-shioda}.

\bibliographystyle{amsplain}
\bibliography{references}

\providecommand{\bysame}{\leavevmode\hbox to3em{\hrulefill}\thinspace}
\providecommand{\MR}{\relax\ifhmode\unskip\space\fi MR }
\providecommand{\MRhref}[2]{%
  \href{http://www.ams.org/mathscinet-getitem?mr=#1}{#2}
}
\providecommand{\href}[2]{#2}
\begin{thebibliography}{10}

\bibitem{jac-of-gen-1-curves}
Sang~Yook An, Seog~Young Kim, David~C. Marshall, Susan~H. Marshall, William~G.
  McCallum, and Alexander~R. Perlis, \emph{Jacobians of genus one curves}, J.
  Number Theor. \textbf{90} (2001), 303--315.

\bibitem{ArtebaniSartiTaki}
Michela Artebani, Alessandra Sarti, and Shingo Taki, \emph{K3 surfaces with
  non-symplectic automorphisms of prime order}, Math. Z. \textbf{268} (2011),
  507--533.

\bibitem{DGPS}
Wolfram Decker, Gert-Martin Greuel, Gerhard Pfister, and Hans Sch\"onemann,
  \emph{{\sc Singular} {4-1-0} --- {A} computer algebra system for polynomial
  computations}, 2016, \url{http://www.singular.uni-kl.de}.

\bibitem{degtyarev-supersingular}
Alex Degtyarev, \emph{Lines in supersingular quartics},  (2016), preprint,
  \verb+arXiv:1604.05836+.

\bibitem{degtyarev16}
\bysame, \emph{Smooth models of singular {$K3$}-surfaces}, Rev. Mat. Iberoam.
  (2016), preprint, \verb+arXiv:1608.06746+.

\bibitem{degtyarev-itenberg-sertoz}
Alex Degtyarev, Ilya Itenberg, and Ali~Sinan Sertöz, \emph{Lines on quartic
  surfaces}, Math. Ann. \textbf{368} (2016), no.~1, 753--809.

\bibitem{GAP4}
The GAP~Group, \emph{{GAP -- Groups, Algorithms, and Programming, Version
  4.8.7}}, 2017, \url{http://www.gap-system.org}.

\bibitem{nikulin-finite-groups}
Viacheslav~V. Nikulin, \emph{Finite groups of automorphisms of {Kählerian}
  {$K3$} surfaces}, Trudy Moskov. Mat. Obshch. \textbf{38} (1979), 75--137
  (Russian).

\bibitem{nikulin}
\bysame, \emph{Integral symmetric bilinear forms and some of their
  applications}, Math. USSR Izvestija \textbf{14} (1980), no.~1, 103--167.

\bibitem{oguiso17}
Keiji Oguiso, \emph{Isomorphic quartic {K3} surfaces in the view of {Cremona}
  and projective transformations}, Taiwanese J. Math. \textbf{21} (2017),
  no.~3, 671--688.

\bibitem{64lines}
S{\l}awomir Rams and Matthias Sch{\"u}tt, \emph{64 lines on smooth quartic
  surfaces}, Math. Ann. \textbf{362} (2015), no.~1, 679--698.

\bibitem{rams-schuett-char2}
\bysame, \emph{At most 64 lines on smooth quartic surfaces (characteristic~2)},
   (2017), preprint, \verb+arXiv:1512.01358v2+.

\bibitem{sagemath}
\relax{The Sage Developers}, \emph{{S}agemath, the {S}age {M}athematics
  {S}oftware {S}ystem ({V}ersion 8.0)}, 2017, \url{http://www.sagemath.org}.

\bibitem{schur}
Friedrich Schur, \emph{Ueber eine besondre {Classe} von {Fl{\"a}chen} vierter
  {Ordnung}}, Math. Ann. \textbf{20} (1882), 254--296.

\bibitem{segre}
Beniamino Segre, \emph{The maximum number of lines lying on a quartic surface},
  Q. J. Math. \textbf{14} (1943), 86--96.

\bibitem{shimada-shioda}
Ichiro Shimada and Tetsuji Shioda, \emph{On a smooth quartic surface containing
  $56$ lines which is isomorphic as a {$K3$} surface to the {Fermat} quartic},
  Manuscripta Math. \textbf{153} (2016), no.~1--2, 279--297.

\bibitem{shioda-inose77}
Tetsuji Shioda and Hiroshi Inose, \emph{On singular {K3} surfaces}, Complex
  analysis and algebraic geometry, Iwanami Shoten, Tokyo, 1977, pp.~119--136.

\bibitem{vanGeemenSarti}
Bert van Geemen and Alessandra Sarti, \emph{Nikulin involutions on {$K3$}
  surfaces}, Math. Z. \textbf{255} (2007), 731--753.

\bibitem{veniani1}
Davide~Cesare Veniani, \emph{The maximum number of lines lying on a {K3}
  quartic surface}, Math. Z. \textbf{285} (2017), no.~3, 1141--1166.

\end{thebibliography}

\end{document}